\documentclass[10pt,twoside]{amsart}
\usepackage[all]{xy}   
        \usepackage {amssymb,latexsym,amsthm,amsmath}
        
\long\def\symbolfootnote[#1]#2{\begingroup%
\def\thefootnote{\fnsymbol{footnote}}\footnote[#1]{#2}\endgroup}

\makeatletter
\def\imod#1{\allowbreak\mkern10mu({\operator@font mod}\,\,#1)}
\makeatother

\newtheorem{theorem}{Theorem}[section]

\newtheorem*{theorem*}{Theorem}
\theoremstyle{definition}

\numberwithin{equation}{section}

\newcommand{\ignore}[1]{}

\newcommand{\mynote}[1]{}

\begin{document}
\setcounter{section}{0}
% document information
\title{Real Elements and Schur Indices of a Group}
\author{Amit Kulsherstha}
\address{IISER Mohali, MGSIPA Complex, Sector-26, Chandigarh 160019 INDIA.}
\email{amitk@iisermohali.ac.in}

\author{Anupam Singh}
\address{IISER Pune, central tower, Sai Trinity building, Pashan circle, Sutarwadi, Pune 411021 INDIA.} 
\email{anupamk18@gmail.com}
\date{}
\thanks{The second named author thanks IMS for invitation to give a talk in its 76th annual meeting held in Surat.}


%%%%%%%%%%%%%%%%%%%%%%%%%%%%

\begin{abstract}
In this article we try to explore the relation between real conjugacy classes and real characters of finite groups at more refined level.
This refinement is in terms of properties of groups such as strong reality and total orthogonality. In this connection we raise several questions
and record several examples which have motivated those questions.
\end{abstract}

\maketitle
%%%%%%%%%%%%%%%%%%%%%%%%%%%%
\section{Introduction}
Let $G$ be a finite group. It is a basic statement in the character theory of finite groups
that the number of irreducible complex characters is same
as the number of conjugacy classes in $G$. Further this statement can be refined to the number of real irreducible
characters being equal to the number of real conjugacy classes. However it is very well known that the real characters
come from two kinds of representations: {\it orthogonal} and {\it symplectic} \cite[Section 13.2 Prop. 39]{se}. Schur himself explored
this topic and now this is very closely related to Schur index computation. In last 10 years the computation
of Schur indices have been almost completed for finite groups of Lie type (for example see \cite{ge}). 

In this article we raise the question to further divide real conjugacy classes in two parts as to match the sizes of partitions
with the number of orthogonal characters and symplectic ones. In section~\ref{conjclasses} we give basic definitions 
and ask several questions concerning relations between conjugacy classes and characters. 
In section~\ref{schurindex} we define Schur indices of a 
representation. In section~\ref{caninv} we define canonical involution on a group algebra and mention known results
about when the involution restricts to the simple components. 
In the following section we provide examples which provide strength
to the questions raised earlier. This is the main objective of this article.
In the section~\ref{plesken} we mention the Lie algebra defined for a group algebra
which makes use of real conjugacy classes. Results mentioned in this article are either calculations using computer algebra
system GAP or have been collected from various sources related to real conjugacy classes. Some of the questions raised in
this article are already known to experts (for example see~\cite{gn}).

%%%%%%%%%%%%%%%%%%%%%%%%%%%%%%%%
\section{Conjugacy Classes vs. Representations for a group}\label{conjclasses}

Let $G$ be a group. An element $g\in G$ is called {\it real} if there exists $t\in G$ such that $tgt^{-1}=g^{-1}$.  
An element $g \in G$ is called {\it involution} if $g^2=1$. An element $g\in G$ is called {\it strongly real} if
it is a product of two involutions in $G$. Further notice that if an element is (strongly) real then all its conjugates are (strongly) real.
Hence reality (i.e. being real) and strong reality are properties of conjugacy classes.
The conjugacy classes of involutions and more generally strongly real elements are obvious examples of real classes.
However converse need not be true.

\begin{example}
Take $G=Q_8$, the finite quaternion group. Then we see that $iji^{-1}=j^{-1}$ hence $j$ is real but it is not strongly real.
\end{example}

A {\it representation} of a group $G$ is a homomorphism $\rho \colon G\rightarrow GL(V)$ where $V$ is a vector space 
over a field $k$. The representation $\rho$ is called {\it irreducible} if $V$ and $\{0\}$ are the only subspaces
$W$ of $V$ satisfying $\rho(G).W \subseteq W$. Let $V^*$ denote the dual vector space of $V$. The {\it dual representation} of
$\rho$ is the representation $\rho^*\colon G\rightarrow GL(V^*)$ given by $\rho^*(g) = ~^t\rho(g^{-1})$, where
$~^t\rho(g^{-1})$ is the transpose of $\rho(g^{-1})$.
To a representation $\rho$ its associated {\it character} $\chi \colon G\rightarrow k$ is
defined by $\chi(g)=trace(\rho(g))$. In this section we only consider complex representations (i.e. $k=\mathbb C$).
It is a classical theorem that the number of conjugacy classes in $G$ is equal to the number of irreducible complex characters
\cite[Section 2.5 Theorem 7]{se}. A theorem due to Brauer \cite[Chapter 23]{jl} asserts that 
the number of real conjugacy classes is same as the number of irreducible real characters (i.e. the complex characters which 
take real values only). 

Let $\rho \colon G\rightarrow GL(V)$ be an irreducible complex representation of $G$ and $\chi$ be the associated character. If $\chi$
takes a complex value then the representation $\rho$ is not isomorphic to its dual $\rho^*$ and the vector space $V$
can not afford a $G$-invariant
non-zero bilinear form. If $\chi$ is real then $\rho\cong \rho^*$ and in this case $V$ admits a
non-zero $G$-invariant bilinear form. This form can be either symmetric or skew-symmetric depending on whether $V$ is defined over
$\mathbb R$ or not (see \cite[Section 13.2]{se}). Equivalently the image of $G$ sits inside $O_n$ in the first case and inside $Sp_{2n}$ in the second case.
To each character $\chi$ one associates {\it Schur indicator} $\nu(\chi)$ which is defined as follows:
$$\nu(\chi)=\frac{1}{|G|}\sum_{g\in G}\chi(g^2).$$
In fact $\nu(\chi)=0,\pm 1$ and $\nu(\chi)=0$ if and only if $\chi$ is not real, $\nu(\chi)=1$ if $\chi$ is
orthogonal and symplectic otherwise \cite[Prop. 39]{se}. Hence the real characters come in two classes: orthogonal
type and symplectic type. This gives a natural division of real characters. We now have following questions:
\begin{question} Let $G$ be a finite group.
\begin{enumerate}
\item Can we naturally divide real conjugacy classes of $G$ in two parts so that the number of one part is same as the
number of orthogonal representations (we are specially interested in groups of Lie type)?
\item Is it true that if a group $G$ has no symplectic character (i.e. all self-dual representations are orthogonal) then 
all real elements are strongly real and vice versa?
\end{enumerate}
\end{question}

\noindent A careful look at the examples in section \ref{example} will suggest that these questions are indeed interesting to ask.

We now restrict our attention to those groups in which all elements are real. Such groups are called {\it real} or {\it ambivalent}.
Tiep and Zalesski classify all real finite quasi-simple groups in \cite{tz}. For real groups all characters are real valued. Two interesting subclasses of the real groups are the following:
\begin{itemize}
\item Subclass in which all elements are strongly real. Groups belonging to this class are called {\it strongly real} groups.
\item Subclass in which all characters are orthogonal. Groups belonging to this class are called {\it totally orthogonal}
or {\it ortho-ambivalent}.
\end{itemize}

An element $g\in G$ is called {\it rational} if $g$ is conjugate to $g^i$ whenever $g$ and $g^i$ generate the same subgroup of $G$.
A group is called {\it rational} if every element of $g$ is rational. 
\begin{theorem}[\cite{se} 13.1-13.2]\label{rational-char} 
A group $G$ is rational if and only if all its characters are $\mathbb Q$-valued.
\end{theorem}
In fact, the number of isomorphism classes of irreducible representations of $G$ over $\mathbb Q$ is same as the number of 
conjugacy classes of cyclic subgroups of $G$.
\begin{example}
Let $G=S_n$, the symmetric group. Then every element of $G$ is a product of two involutions and hence strongly real. It is
also totally orthogonal. All its character are integer valued. In fact, it is a rational group.
\end{example}

\begin{example}
The group $Q_8$ has four $1$-dimensional representation which are orthogonal and one symplectic representation. This group
is neither strongly real nor totally orthogonal.
\end{example}

Let us denote the class of finite groups which are real by $\mathcal R$, the class of real groups which have their
Sylow $2$-subgroup Abelian by $\mathcal S$, the class of strongly real groups by $\mathcal{SR}$ and the class of 
totally orthogonal groups by $\mathcal{TO}$. Then using the results of~\cite{wg} and~\cite{ar} we have the following:
$$\mathcal S \subset \mathcal{SR} \subset \mathcal R  \ \ \ \ \text{and}\ \ \ \  
\mathcal S \subset \mathcal{TO} \subset \mathcal R.
$$
In particular they prove that if $G$ is real then it is generated by its $2$-elements and if $G$ is totally orthogonal
then it is generated by involutions. In view of above relations we ask following question:
\begin{question}
Find the class $\mathcal{SR}\cap \mathcal{TO}$. 
\end{question}

We know that the class $\mathcal{SR}\cap \mathcal{TO}$ contains $\mathcal S$. The containment though, is far
from being equality. The dihedral group $D_4$ of order $8$ is strongly real as well as totally orthogonal but
its Sylow $2$-subgroup, which is $D_4$ itself, is not Abelian. 

It is therefore natural to ask which strongly real groups are totally orthogonal and vice
versa. We have made many calculations using GAP and it seems the class $\mathcal{SR}\cap \mathcal{TO}$ is very
close to the class $\mathcal{TO}$,
though $\mathcal{SR}$ and $\mathcal{TO}$ are not identical. There is a group of order $32$ which is strongly real 
but not totally orthogonal. This group $G$ has the following properties: 

\begin{enumerate}
\item It is not simple. It has normal subgroups of all plausible orders - $2$, $4$, $8$ and $16$.
\item Exponent of $G$ is $4$.
\item Derived subgroup of $G$ of order $2$.
\item The group $G$ is a semidirect product of $C_2 \times Q_8$ and $C_2$. Here $C_2$ denotes the
cyclic group of order $2$ and $Q_8$ denotes the quaternion group of order $8$.
\item This group has one $4$-dimensional character and sixteen $1$-dimensional characters. The
$4$-dimensional character assumes non-zero value only on one non-identity conjugacy class.
\end{enumerate}

It is worth noting that this is the only group of order $63$ or smaller which is strongly
real and not totally orthogonal. All totally orthogonal groups till this order are
strongly real.

%%%%%%%%%%%%%%%%%%%%%%%%%%%%%%%%%%%
\section{Schur Indices of groups}\label{schurindex}

Let $G$ be a finite group. Let $k$ be a field such that $char(k)$ does not divide $|G|$. The {\it group algebra} of 
$G$ over $k$ is $\displaystyle kG=\left\{\sum_{g\in G} \alpha_gg : \alpha_g\in k\right\}$ with operations defined as follows:
\begin{eqnarray*}
\sum_{g\in G} \alpha_gg + \sum_{g\in G} \beta_gg &=& \sum_{g\in G} (\alpha_g+\beta_g)g \\
\alpha \left(\sum_{g\in G} \alpha_gg\right) &=& \sum_{g\in G} \alpha\alpha_gg \\
\left(\sum_{g\in G} \alpha_gg\right).\left(\sum_{h\in G} \beta_hh\right) &=& \sum_{t\in G} \left(\sum_{g}\alpha_g\beta_{g^{-1}t}\right)t
\end{eqnarray*}
It is a classical result of Maschke that $kG$ is a semisimple algebra. Hence by using Artin-Wedderburn theorem one can write it 
as a product of simple algebras over division algebras, i.e.,
$$kG\cong M_{n_1}(D_1)\times\ldots\times M_{n_r}(D_r)$$
where $D_i$'s are division algebras over $k$ with center, say $L_i$, a finite field extension of $k$. 
We know that $D_i$ over $L_i$ is of square dimension, say $m_i^2$. Further each simple component
$M_{n_i}(D_i)$ corresponds to an irreducible representation of $G$ over $k$.
The number $m_i$ is called the {\it Schur index} of the corresponding representation. One can write $1=e_1+\cdots +e_r$ using 
above decomposition where $e_i$'s are idempotents. In fact,
$$e_i = \frac{\chi_i(1)}{|G|} \sum_{g\in G} \chi_i(g^{-1})g.$$
where $\chi_i$ is the character of the representation corresponding to $M_{n_i}(D_i)$.
For this reason often the Schur index is denoted as $m_k(\chi_i)$. 
The following are important questions in the subject:
\begin{question}
\begin{enumerate}
\item Find out Schur indices $m_k(\chi_i)$ and $L_i$ for the representations of a group $G$. 
\item Determine the division algebras $D_i$ appearing in the decomposition. 
\end{enumerate}
\end{question} 

\noindent The problem of determining Schur indices for groups of Lie type has been studied extensively in the litrature notably 
by Ohmori \cite{oh1, oh2}, Gow \cite{go1,go2}, Turull (\cite{tu}) and Geck(\cite{ge}). 
Geck also gives a table of the all known results on page 21 in \cite{ge}.
However answer to the second question is much more difficult, e.g., Turull does it for $SL_n(q)$ in \cite{tu}.

The group algebra is well studied over field $k=\mathbb C,\mathbb R, \mathbb Q$ or $\mathbb F_q$. For example, $\mathbb CG\cong 
M_{n_1}(\mathbb C)\times\ldots\times M_{n_r}(\mathbb C)$, as there is only one finite dimensional division algebra over $\mathbb C$ which
is $\mathbb C$ itself (so is over $\mathbb F_q$). Moreover the simple components in this decomposition correspond to a finite dimensional representation
of $G$ over $\mathbb C$.
However in the case of $\mathbb R$ we know that the finite dimensional division algebras over $\mathbb R$ are 
$\mathbb R,\mathbb C$ or $\mathbb H$ hence the corresponding Schur index is $1,2$ or $1$ respectively.
In this case, $$\mathbb RG\cong M_{n_1}(\mathbb R)\times\ldots\times M_{n_l}(\mathbb R)\times M_{n_{l+1}}(\mathbb H)\times
\ldots M_{n_{l+s}}(\mathbb H)\times M_{n_{l+s+1}}(\mathbb C)\times \ldots M_{n_{l+s+p}}(\mathbb C).$$
In particular we see that the irreducible representations of $G$ are of three kinds which are called orthogonal,
symplectic and unitary as they can afford a symmetric bilinear form, an alternating form or a hermetian form respectively. 
These corrspond to the Schur indicator $\nu(\chi)$ (defined in section~\ref{conjclasses}) being $1, -1$ or $0$ respectively. 
Hence the question of calculating Schur indices over $\mathbb R$ is related to determining the types of
representations: orthogonal, symplectic or unitary.

%%%%%%%%%%%%%%%%%%%%%%%%%%%%%%%%%%%%%%%%%%%%%%%
\section{Canonical Involution on the Group Algebra}\label{caninv}

One can define an involution $\sigma$ on $kG$ as follows: 
$$\sigma \left(\sum\alpha_gg \right)\mapsto \sum\alpha_gg^{-1}.$$
This involution is called the {\it canonical involution}. 
We can define a symmetric bilinear form $T\colon kG \times kG\rightarrow k$ by $T(x,y)=tr(l_{x\sigma(y)})$ where $l_x$ is the
left multiplication operator on $kG$ and $tr$ denotes its trace.
We note that $tr(l_e)=n$ and $tr(l_g)=0$ for $g\neq e$.
Hence the elements of group $G$ form an orthogonal basis and the form $T \simeq n<1,1,\ldots,1>$. 
The following result is proved in \cite{sc} (chapter 8 section 13): If the form $n<1,1,\ldots,1>$ is anisotropic over $k$ 
(for example $k=\mathbb R$ or $\mathbb Q$) the involution $\sigma$ restricts to each simple component of $kG$.
In fact, in the case (ref. \cite{bo} theorem 2) $G$ is real the involution $\sigma$ restricts to each simple component of $kG$. 

In general if $kG\cong A_1\times A_2\times\cdots\times A_r$ and $k=\mathbb C$ then either $\sigma$ restricts to a component
$A_i$ or it is a switch involution on $A_i\times A_j$ where $A_i\cong A_j$.
However when $k=\mathbb R$ the involtuion $\sigma$ restricts to each component $A_i$, say $\sigma_i$, and is of either 
first type or second type. It is of the first type when $A_i$ is $\cong M_n(\mathbb R)$ or $M_m(\mathbb H)$ and of the second 
type when the component is isomorphic to $M_l(\mathbb C)$. Moreover when we tensor this with $\mathbb C$ the first type
gives the component over $\mathbb C$ the one to which $\sigma$ restricts and the second type over $\mathbb R$ gives the
one which correspond to the switch involution over $\mathbb C$. 

Let us assume now that the canonical involution $\sigma$ restricts to all components of $kG$, i.e., 
$(kG,\sigma)\cong \prod_i (A_i,\sigma_i)$ where $A_i$s are simple algebra over $k$ with involution. For example, this happens
when $k=\mathbb R$ or when $G$ is real. Algebras with involution
 $(A,\sigma)$ have been studied in the litrature (see \cite{kmrt}) extensively for its connection with algebraic groups. They
are of two kinds: {\it involution of the first kind} is one which restricts to the center of $A$ as identity and the {\it involution of 
second kind} restrict to the center of $A$ as order $2$ element. Further the involution of the first kind are of two types called
{\it orthogonal type} and {\it symplectic type} the second kind is also called {\it unitary type}.

A group is called {\it ortho-ambivalent} with respect to a field $k$ if the canonical involution $\sigma$
restricts to all of its simple components as orthogonal involution (of first kind).
The following results are proved in the thesis of Zahinda (ref. \cite{za} Chapter 2):
An ortho-ambivalent group is necessarily ambivalent and in fact it is totally orthogonal (see Proposition 2.4.2 in \cite{za}). 
The notion of ortho-ambivalence over $k$ is equivalent to ortho-ambivalence over $\mathbb C$. Further the question that which $2$ groups are ortho-ambivalent
is analysed.

%%%%%%%%%%%%%%%%%%%%%%%%%%%%%%%%%%%%%%%%%%%%%%%%%%%%
\section{Some Examples}\label{example}
Here we write down some examples and some GAP calculations we did.

\subsection{Symmetric and Alternating Groups}

Conjugacy classes in $S_n$ are in one-one correspondance with partitions of $n$. Every conjugacy class in $S_n$ is strongly real and hence real.
All characters of $S_n$ are real and moreover orthogonal.

Let $g\in A_n$. Then if $\mathcal Z_{S_n}(g)\subset A_n$ then $g^{S_n}=g^{A_n} \cup (xgx^{-1})^{A_n}$ where $x$ is an 
odd permutation. In case $\mathcal Z_{S_n}(g)\not\subset A_n$ then $g^{S_n}=g^{A_n}$.
In \cite{pa}, Parkinson classified real elements in $A_n$.
Let $n=n_1+n_2+\cdots+n_r$ be a partition and $C$ be a conjugacy class corresponding to that in $S_n$ contained in $A_n$. 
Then, $C$ is real in $A_n$ if and only if 
\begin{enumerate}
\item each $n_i$ distinct,
\item each $n_i$ odd and
\item $\frac{1}{2}(n-r)$ is odd.
\end{enumerate}
And hence the number of conjugacy classes in $A_n$ is equal to the number of real even partitions $+$ twice the number 
of non-real even partitions.
In fact (see \cite{br,pa}), $A_n$ is ambivalent if and only if $n=1, 2, 5, 6, 10, 14$.

For example in the case of $n=7$, the partitions are given by $1^7, 1^62, 1^43, 1^32^2, 1^223, 13^2, 2^23, \\ 34, 25,  124,  16, 7,
1^25, 1^34, 12^3$. Out of which $1^7, 1^43, 1^32^2, 13^2, 2^23, 124,  7, 1^25$ correspond to elements in $A_7$. By above
criteria the only non-real class corresponds to the partition $7$. Hence there are $7$ real conjugacy classes in $A_7$ out of total $7+2.1=9$
conjugacy classes.
Using GAP we verified following statements about $A_7$.
\begin{enumerate}
\item All but one conjugacy classes are real.
\item All real conjugacy classes are strongly real.
\item All real characters are orthogonal.
\end{enumerate}

We summarise some GAP claculations below for $A_n$:
\vskip2mm
\begin{center}
\begin{tabular}{|c|c|c|c|c|c|c|}\hline
n& total &real classes & st real & orthogonal & symplectic & unitary\\
 & classes &             &          & characters & characters & characters\\
 \hline
 5 & 5 & 5 & 5   &  5   & 0 & 0 \\ \hline
 6 & 7 & 7 & 7   &  7   &  0 & 0 \\ \hline
 7 & 9 & 7 & 7   &  7   &  0 & 2 \\ \hline
 8 & 14 & 10 & 10 & 10 & 0 & 4 \\ \hline
 9 & 18 & 16 & 16 & 16 & 0 & 2  \\\hline
 10 & 24 & 24 & 24 & 24 &0 & 0\\ \hline
 14 & 72 & 72 & 72 & 72 & 0 & 0 \\ \hline
\end{tabular}
\end{center}
\vskip3mm

In~\cite{su} section 3, Suleiman proved that in alternating groups an element is real if and only if it is strongly real.
Moreover in $A_n$ every real character is orthogonal, i.e., the Schur index of $A_n$ is $1$.
This follows from a work of Schur which is quoted in a paper of Turull (see \cite{tu2}, Theorem 1.1).
This result of Schur says that the Schur index of every ireducible representation of $A_n$ (for {\it each} $n$) is $1$. 
Thus $A_n$ has no symplectic characters. Hence for $A_n$,\\

$|$strongly real classes$|$ = $|$real classes$|$ =  $|$real characters$|$ = $|$orthogonal characters$|$. \\

Here the first equality is from Suleiman, second is obvious and third one is the Schur index computation by Schur himself.

\vskip3mm
Now one can ask similar questions for the covers of these groups namely $\tilde S_n$ and $\tilde A_n$.
We summarise some GAP claculations below for $\tilde{A_n}$:
\vskip2mm
\begin{center}
\begin{tabular}{|c|c|c|c|c|c|c|}\hline
n& total &real classes & st real & orthogonal & symplectic & unitary\\
& classes &             &          & characters & characters & characters\\
 \hline
4 & 7 & 3 & 2  &  2   & 1 & 4 \\ \hline 
5 & 9 & 9 & 2  &  5   & 4 & 0 \\ \hline 
\end{tabular}
\end{center}
\vskip3mm

We summarise some GAP claculations below for $\tilde{S_n}$:
\vskip2mm
\begin{center}
\begin{tabular}{|c|c|c|c|c|c|c|}\hline
n& total &real classes & st real & orthogonal & symplectic & unitary\\
& classes &             &          & characters & characters & characters\\
 \hline
4 & 8 & 6 & 6  &  6   & 0 & 2 \\ \hline 
5 & 12 & 8 & 6  &  7   & 1 & 4 \\ \hline 
\end{tabular}
\end{center}
\vskip3mm

Here we see an example of group, say the Schur cover of $A_5$ for which there are only $2$ strongly real classes and all $9$ real classes, 
while it has $5$ orthogonal characters and remaining $4$ symplectic ones.

\subsection{$GL_n(q)$}
The group $GL_n(q)$ has the property that all real elements are strongly real (ref. \cite{wo}). 
It doesn't have irreducible symplectic representations, i.e., all self-dual irreducible representations 
are orthogonal (\cite{dp} Theorem 4). In \cite{gs1} and \cite{gs2} the precise number of the real elements are calculated.
In \cite{ma}, Macdonald gives an easy way to enumerate conjugacy classes.
\begin{theorem}[Macdonald]
Conjugacy classes in $GL_n(q)$ are in one-one correspondance with a sequence of polynomials $u=(u_1,u_2,\ldots)$ satisfying:
\begin{enumerate}
\item a partition of $n$, $\nu=1^{n_1}2^{n_2}\cdots$, i.e., $|\nu|=\sum_i in_i=n$,
 \item $u_i(t)=a_{n_i}t^{n_i}+\cdots+a_1t+1 \in \mathbb F_q[t]$ for all $i$ with $a_{n_i}\neq 0$.
\end{enumerate}
\end{theorem}
\noindent Hence the number of conjugacy classes in $GL_n(q)$ is 
$$\sum_{\{\nu : |\nu|=n\}} c_{\nu}= \sum_{\{\nu : |\nu|=n\}}\prod_{n_i>0}(q^{n_i}-q^{n_i-1}).$$
\begin{theorem}[\cite{gs1}]
Real conjugacy classes in $GL_n(q)$ are in one-one correspondance with a sequence of polynomials $u=(u_1,u_2,\ldots)$ satisfying
\begin{enumerate}
\item a partition of $n$, $\nu=1^{n_1}2^{n_2}\cdots$, i.e., $|\nu|=\sum_i in_i=n$,
 \item $u_i(t)=a_{n_i}t^{n_i}+\cdots+a_1t+1 \in \mathbb F_q[t]$ for all $i$ with $a_{n_i}\neq 0$.
\item $u_i(t)$ self-reciprocal.
\end{enumerate}
\end{theorem}
\noindent Hence the number of real conjugacy classes in $GL_n(q)$ is $$\sum_{\{\nu : |\nu|=n\}}\prod_{n_i>0}n_{q,n_i}$$ where 
$n_{q,n_i}$ is the number of polynomials $u_i(t)$ of above kind of degree $n_i$ over field $\mathbb F_q$.

Hence in $GL_n(q)$ the number of strongly real elements  is same as the number of orthogonal characters.

\subsection{$SL_2(q)$}
In this case if $q$ is even all $q+1$ classes are real as well as strongly real.
If $q$ is odd, there are exactly $2$ strongly real classes. In the case $q\equiv 1 \imod 4$ all $q+4$ classes are real and
if $q\equiv 3 \imod 4$ only $q$ out of $q+4$ are real (in fact, exactly unipotent ones are not real).
Hence we can say that the groups $SL_2(q)$ are ambivalent if and only if $q$ is a sum of (at most) two squares.

In the case $q$ is even all characters are orthogonal. However if $q$ is odd there is always a symplectic character. One can refer 
to the calculations of Schur indices in~\cite{sh} Theorem 3.1. Hence one can conclude that the group $SL_2(q)$ 
is ortho-ambivalent if and only if $q$ is even.

\subsection{$SL_n(q)$}
The real and strongly real classes are calculated in \cite{gs1}. Turull calculated Schur indices of characters of $SL_n(q)$ in~\cite{tu} 
 over $\mathbb Q$ and also determined the division algebras appearing in the decomposition of group algebra.  
The following is known from the work of Ohmori, Gow, Zelevinsky, Turull, Geck etc.
For the group $SL_n(q)$ if $n$ is odd or $n\equiv 0 \imod 4$ or $|n|_2>|p-1|_2$ then all real characters are orthogonal. 
In the case $2\leq |n|_2\leq |p-1|_2$ there are symplectic representations.

\subsection{Orthogonal Group}
In~\cite{wo}, Wonenburger proved that every element of orthogonal group is a product of two involutions. Hence these groups 
are strongly real. In~\cite{go3}, Gow proved (Theorem 1) that all characters of $O_n(q)$ are orthogonal.

%%%%%%%%%%%%%%%%%%%%%%%%%%%%%%%%%%%%%%%%%%%%%%
\section{The Lie Algebra $\mathcal L(G)$}\label{plesken}

Since $kG$ is an associative algebra we can define $[x,y]=xy-yx$ which makes it a Lie algebra. The subspace
$\mathcal L(G)$ generated by $\{\Hat g=g-g^{-1} \mid g\in G\}$ is Lie subalgebra. This Lie algebra associated to
a finite group is studied in~\cite{ct} and called Plesken Lie algebra. They prove,
\begin{theorem}
The Lie algebra $\mathcal L(G)$ admits the decomposition:

$$\mathcal L(G)\cong \bigoplus_{\chi\in \mathbb R} o(\chi(1)) \oplus\bigoplus_{\chi\in \mathbb H} sp(\chi(1)) \oplus
\bigoplus'_{\chi\in \mathbb C}gl(\chi(1))$$
where the sums are over different kind of irreducible characters (with obivious meaning) and the last sum is $'$ed meaning
we have to take only one copy of $gl(\chi(1))$ for $\chi$ and $\chi^{-1}$.
\end{theorem}
They also prove that the Lie algebra $\mathcal L(G)$ is semisimple if and only if  $G$ has no complex characters and every
character of degree $2$ is of symplectic tpe, i.e., $G$ is ambivalent and every non-linear character is symplectic.
They also classify when $\mathcal L(G)$ is simple.

%%%%%%%%%%%%%%%%%%%%%%%%%%%%%%%%%%%%%%%%%%%%%%%
\section{Group Algebra and Real Characters}

We have $\mathbb RG$ a group algebra with involution $\sigma$ which restricts to each simple component of it. On one 
hand we see that $\mathcal Z(\mathbb RG) \cong \oplus_g \mathbb R c_g$ where sum on the right hand side is over conjugacy
classes and $c_g=\sum_{t\in G} g^t$ and on other hand we have $\mathcal Z(\mathbb RG)\cong \oplus \mathcal Z(M_n(D)) \cong 
\bigoplus_{\chi\in\mathbb R}\mathbb R\oplus \bigoplus_{\chi\in\mathbb C} \mathbb C \oplus \bigoplus_{\chi\in\mathbb H} 
\mathbb R$. We know that center of a semisimple algebra is an \'{e}tale algebra and $\sigma$ restricts to it. In fact, in 
the first situation we have $\mathcal Z(\mathbb RG) \cong \oplus_{g\sim g^{-1}} \mathbb R c_g \oplus (\mathbb Rc_g \oplus 
\mathbb Rc_{g^{-1}})$ where on the first component $\sigma$ restricts as trivial involution and on the second component
it becomes as a switch involution (hence of the second kind). In the second isomorphism we know that $\sigma$ restricts to
the trivial map on the $\mathbb R$ and $\mathbb H$ components and is of the second kind on $\mathbb C$ components. Hence
counting the components where $\sigma$ restricts as first kind gives us the number of real conjugacy classes is same as
the number of real plus symplectic representations. 

However applying the same trick to $\mathbb QG$ doesn't give the corresponding result regarding rational representations.
As it will also count the odd degree field extensions of $\mathbb Q$. It will be interesting to find such a proof.

%%%%%%%%%%%%%%%%%%%%%%%%%

\end{document}